\documentclass[11pt]{article}
\usepackage{graphicx}
\usepackage{subfigure}
\usepackage{color}
\usepackage{amsmath,amssymb}
\def\la{\Big\langle}
\def\ra{\Big\rangle}

\def\forall{\hbox{for all}~~}
\def\bfe{{\bf e}}
\def\bfk{{\bf k}}
\def\L{{\bf L}}

\def\E{{\cal E}}

\def\bfv{{\bf v}}
\def\bfa{{\bf a}}

\def\bfc{{\bf c}}
\def\bfh{{\bf h}}

\def\bfw{{\bf w}}

\def\ds{\displaystyle}

\def\bfn{{\bf n}}
\def\ve{\varepsilon}

\def\D{{\cal D}}
\def\N{{\cal N}}

\def\R{{\mathbb R}}

\def\v{\vskip 1em}
\def\O{{\cal O}}

\def\C{{\cal C}}

\def\ov{\overline}

\def\bel{\begin{equation}\label}
\def\eeq{\end{equation}}
\def\sqr#1#2{\vbox{\hrule height .#2pt
\hbox{\vrule width .#2pt height #1pt \kern #1pt
\vrule width .#2pt}\hrule height .#2pt }}
\def\square{\sqr74}
\def\endproof{\hphantom{MM}\hfill\llap{$\square$}\goodbreak}
\def\bega{\begin{array}}
\def\enda{\end{array}}
\def\begi{\begin{itemize}}
\def\endi{\end{itemize}}

\begin{document}
\title{\bf Well-posedness of a Model for the Growth of 
Tree Stems and Vines}

\author{Alberto Bressan and Michele Palladino\\ \, \\
Department of Mathematics, Penn State University.\\
University Park, PA~16802, USA.\\
\\
e-mails:~axb62@psu.edu,
~mup26@psu.edu}
\maketitle
\begin{abstract}

The paper studies a PDE model for the growth 
of a tree stem or a vine, having the form of a differential inclusion 
with state constraints.  The equations
describe the elongation due to cell growth, and the response to gravity 
and to external obstacles. 

The main theorem shows that the evolution problem is well posed, 
until a specific ``breakdown configuration" is reached.
A formula is proved, characterizing the reaction produced by unilateral 
constraints.   At a.e.~time $t$, this is determined by the minimization of
an elastic energy functional under suitable constraints.
\end{abstract}

\section{Introduction}
\label{s:1}
\setcounter{equation}{0}
We consider a PDE model, recently introduced in \cite{BPS},
describing the growth of a plant stem or a vine.

The position of the stem
at time $t$ is described by a curve $\gamma(t,\cdot)$.
For $s\in [0,t]$, we think of $\gamma(t,s)$ as the position
at time $t$ of the cell born at time $s$.   
The model takes into account:
\begi
\item[{\bf (1)}] the linear elongation,
\item[{\bf (2)}]the upward bending, as a response to gravity,
\item[{\bf (3)}] an additional bending, in case of a vine clinging
to branches of other plants,
\item[{\bf (4)}] the reaction produced by obstacles, such as rocks, 
trunks or branches of 
other trees.
\endi
For simplicity, we rescale time and assume that the map $s\mapsto
\gamma(t,s)$ parameterizes the curve by arc-length.
Without loss of generality, one can assume that $\gamma(t,0)=0\in\R^3$,
so that
\bel{g}\gamma(t,s)~=~\int_0^s \bfk(t,\sigma)\, d\sigma\,,
\qquad\qquad \bfk(t,s)~\doteq~\gamma_s(t,s).
\eeq

\begin{figure}[ht]
\centerline{\hbox{\includegraphics[width=8cm]{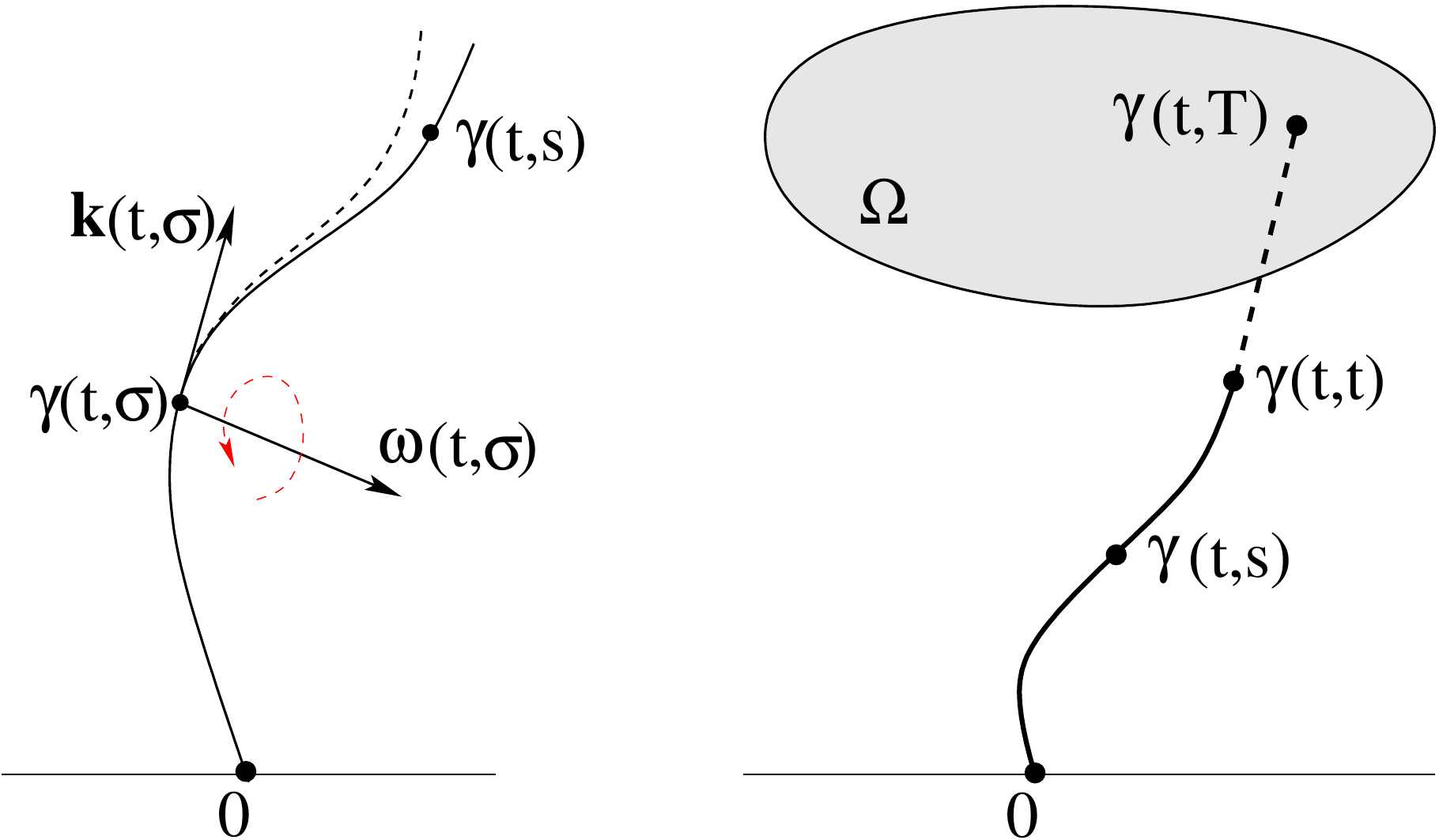}}}
\caption{\small Left: at any point $\gamma(t,\sigma)$ along the stem, 
an infinitesimal change in curvature is produced as a response to gravity (or stems of other plants). The angular velocity is given by the vector $\omega(\sigma)$. This affects the position of all higher points along the stem.
Right: At a given time $t$, the curve $\gamma(t,\cdot)$ is parameterized by $s\in [0,t]$.
It is convenient to 
prolong this curve by adding a segment of length $T-t$ at its tip
(dotted line, possibly entering inside the obstacle).
This yields an evolution equation on a fixed functional space $H^2([0,T];
\,\R^3)$.}
\label{f:sg102}
\end{figure}
\v
The change in the position of points on  the stem
is described by 
\bel{F1}
\ds\gamma_t
(t,s)~=~\int_0^s  \omega(t,\sigma)\times\bigl(\gamma
(t,s)-\gamma
(t,\sigma)\bigr)
\, d\sigma~ \doteq~ F(t,s)\,.\eeq
Here $\omega$ represents an angular velocity
(see Fig.~\ref{f:sg102}).   According to (\ref{F1}), portions of the 
stem can slightly change their curvature in time, as a response to gravity 
or (in the case of vines) to branches of other plants.
Notice that the infinitesimal change in curvature at the point $\gamma(t,\sigma)$
affects all  the upper portion of the stem, 
i.e.~all points $\gamma(t,s)$ with 
$s\in [\sigma, t]$. 
In our model, 
$$\omega(t,s)~=~\Psi\bigl(t,s,\gamma(t,s), \gamma_s(t,s)\bigr)$$
depends on the position and on 
the orientation of the stem, at a given point.
For example, to model the bending of the stem in the upward direction
(as a response to gravity),  one can take 
\bel{PSD}
\Psi(t,s, \gamma
, \bfk)~\doteq~e^{-\beta(t-s)}\bfk\times \bfe_3\,.\eeq
Here $\beta>0$ is a stiffness constant,  while
$\bfe_3\in\R^3$ is a unit vector, oriented in the upward vertical direction. Notice that $t-s$ is the age at time $t$ of the cell
born at time $s$. The factor $e^{-\beta(t-s)}$ accounts for the fact that 
older portions of the stem become more stiff, hence their curvature
changes more slowly.

In addition, we consider an obstacle $\Omega\subset\R^3$, whose presence 
imposes the unilateral constraint
\bel{UC}\gamma(t,s)~\notin~\Omega\qquad\qquad \forall s\in [0,t]\,.\eeq
As  in \cite{BPS}, the evolution of the stem can be described by an
equation of the form
\bel{DI}
\gamma_t(t,s)~= ~F(t,s) + \bfv(t,s),\qquad\qquad \bfv(t,\cdot)\in \Gamma(t) ,\eeq
where $\Gamma(t) $ is a cone of admissible velocities 
determined by the constraint reaction.

Under natural assumptions, 
the main theorem in \cite{BPS} provides the existence of a solution to
(\ref{DI}).  This solution is  defined 
up to the first time where a ``breakdown configuration"
is reached, characterized
at (\ref{bad1})-(\ref{bad2}). Examples are shown in Fig.~\ref{f:sg101}.   
The theorem is proved by writing  the evolution equation for $\gamma$
 in the form of a differential inclusion with closed convex right hand side,
in the functional space $H^2([0,T];\,\R^3)$.
The uniqueness of these solutions, however, had remained 
an open question.  

We remark that most of the literature on differential inclusions
with constraints has been concerned with problems of the form
$${d\over dt} x(t)~\in ~F(x(t)) - N_S(x(t)),\qquad\qquad x(t)\in S,$$
where $N_S(x)$ is the outer normal cone to the set $S$ at the point 
$x$.  When the set $S=S(t)$ is allowed to depend on time, this 
is called a  ``perturbed sweeping process", see \cite{CG, CMM, M, RS}.  
In this setting, the Cauchy problem  usually has a unique solution, 
continuously depending on the initial data.

\begin{figure}[ht]
\centerline{\hbox{\includegraphics[width=15cm]{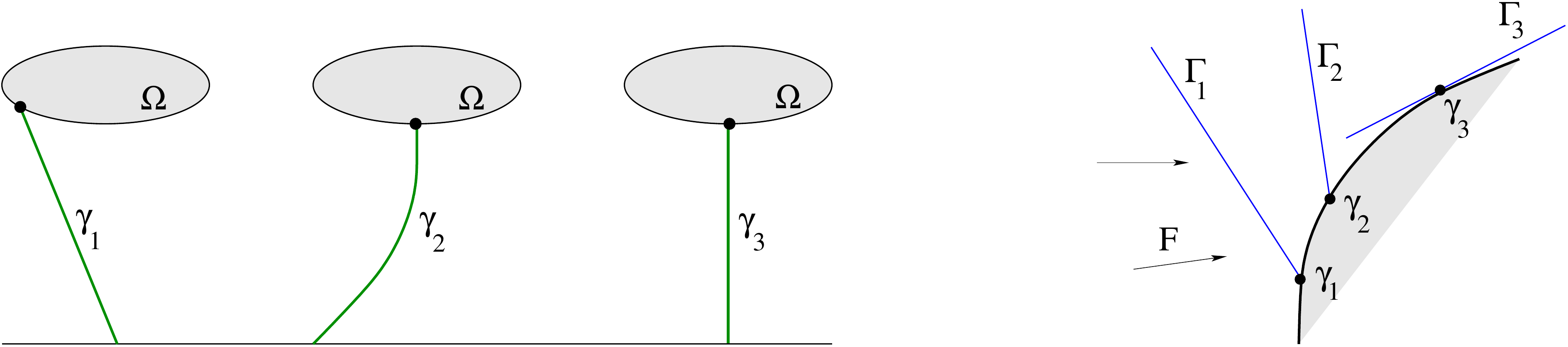}}}
\caption{\small Left: three configurations of the stem, relative to the obstacle.   Right: in an abstract space, the first two 
configurations are represented by points $\gamma_1,\gamma_2$ on the boundary 
of the admissible set $S$ where the corresponding cones $\Gamma_1,\Gamma_2$ are transversal.
On the other hand, $\gamma_3$ is a ``breakdown configuration",
satisfying all assumptions (\ref{bad1})-(\ref{bad2}).
Its corresponding cone $\Gamma_3$ is tangent to the boundary
of the set $S$.  Here the shaded region is the complement of $S$.}
\label{f:sg103}
\end{figure}

On the contrary, in the present case the cone $\Gamma$ of admissible velocities in (\ref{DI})
bears no relation to the normal cone.   In fact, as the stem reaches
a ``breakdown
configuration" illustrated in Fig.~\ref{f:sg103}, the cone $\Gamma$ becomes
tangent to the boundary of the admissible set $S$.
For this reason, the well-posedness of the Cauchy problem for (\ref{DI})
is a delicate issue.    

The aim of the present paper is twofold:
\begi
\item[(i)] Prove the  uniqueness and continuous dependence on initial data,
for solutions to (\ref{DI}).
\item[(ii)] Provide a characterization of 
the velocity  $\bfv(t,\cdot)\in \Gamma(t)$ 
in (\ref{DI}) determined by the obstacle reaction.
\endi
Following \cite{BPS}, a solution $t\mapsto \gamma(t,\cdot)$
is regarded as a map taking values 
in the Hilbert space $H^2([0,T];\,\R^3)$.
Unfortunately, a study of the $H^2$ distance between two solutions does 
not lead to any useful estimate.
In the present paper, 
the distance between two solutions $\gamma_1(t,\cdot)$, $\gamma_2(t,\cdot)$
will be estimated by constructing a family of rotations, transforming 
a unit tangent vector $\bfk_1(t,s)$ to $\gamma_1$ into the corresponding 
 tangent vector  $\bfk_2(t,s)$ to $\gamma_2$, for every $s\in [0,t]$.
By estimating how the norm of these rotation vectors grows in time, 
we shall provide a bound on the distance between the two solutions 
$\gamma_1,\gamma_2$ 
for all times $t$.

Next, by further developing the analysis in \cite{BPS} we will show
that, for a.e.~time $t$, the vector $\bfv(t,\cdot)$ is uniquely determined
by the solution of a variational problem. Indeed, $\bfv$ can be recovered by
the formula (\ref{Min4}), where $\bar\omega(\cdot)$ is
 the minimizer of an elastic 
deformation energy, subject to the unilateral constraints posed by the obstacle $\Omega$.

The remainder of the paper is organized as follows.
In Section~2 we review the model equations and all the main 
definitions and assumptions.
We then recall the existence theorem proved in \cite{BPS}, and state
the main results of the paper; namely, the uniqueness and 
characterization of solutions, stated in Theorems~2 and 3, respectively.
Section~3 contains some preliminary lemmas, on the existence of rotation
vectors transforming one curve into another one.
The uniqueness of solutions is proved in Section~4, while the representation
formula (\ref{Min4}) is proved in Section~5.

For the general theory of optimal control, also in the presence 
of state constraints, we refer to \cite{Cesari, V}. 
A description of plant development from a biological 
point of view can be found in \cite{LD}.

\section{Statement of the main results}
\label{s:2}
\setcounter{equation}{0}
We start with a brief review of the model considered in \cite{BPS}.

At each time $t$, the position of the stem is described by a map
$s\mapsto \gamma(t,s)$ from $[0,t]$ into $\R^3$.    
Clearly, the domain of this map grows with
time.   
It is convenient to reformulate the model as an evolution problem 
on a functional space independent of $t$.   For this purpose, we fix $T>t_0$
and consider the Hilbert-Sobolev space $H^2([0,T];\, \R^3)$.  
Any function $\gamma(t,\cdot)\in H^2([0,t];\, \R^3)$ will be canonically extended to 
$H^2([0,T];\, \R^3)$ by setting (see Fig.~\ref{f:sg102}, right)
\bel{Pex}
\gamma(t,s)~\doteq~\gamma(t,t) + (s-t)\gamma_s(t,t)
\qquad\qquad \hbox{for}~~s\in [t,T]\,.\eeq
Notice that the above extension is well defined because 
$\gamma(t, \cdot)$ and $\gamma_s(t,\cdot)$
are continuous functions.  Throughout the following, we
shall study
functions defined on a domain of the form
\bel{DT}\D_T~\doteq~\bigl\{(t,s)\,;~~0\leq s\leq t,~~
t\in [t_0,T]\bigr\}, \eeq
and extended to the rectangle $[t_0,T]\times [0,T]$ as in (\ref{Pex}).
In particular, the partial derivative $\gamma_s(t,s)$ will be constant 
for $s\in [t,T]$.

Adopting the notation $a\wedge b\doteq \min\{a,b\}$,
we consider an  evolution problem on the space $H^2([0,T];\R^3)$,
having the  form
\bel{E1}
\gamma_t(t,s)~=~\int_0^{s\wedge t} \Psi\Big(t,\sigma, \gamma(t,\sigma), 
\gamma_s(t,\sigma)\Big)\times 
\Big(\gamma(t,s)-\gamma(t,\sigma)\Big)
d\sigma
+\bfv(t,s).\eeq
Here $s\in [0,T]$,  $\Psi: \R\times \R\times \R^3\times\R^3\mapsto \R^3$ is a smooth function, and $\bfv(t,\cdot)$ is an admissible velocity field produced by the constraint reaction.   More precisely, 
let $\Omega\subset\R^3$ be an open set with $\C^2$ boundary.
Given the configuration
$\gamma(t,\cdot)$ of the stem at time $t$, let
\bel{chit}
\chi(t)~\doteq~\Big\{ s\in [0,t]\,;~~
\gamma
(t,s)\in ~\partial\Omega\Big\}
\eeq
be the set where the stem touches the obstacle.
For $s\in \chi(t)$, let $\bfn(t,s)$ be the unit outer normal to the 
boundary $\partial\Omega$ at the point $\gamma
(t,s)$.  
The {\bf cone of admissible velocities} produced by the 
obstacle reaction is defined 
to be the set of velocity fields
\bel{RC} \bega{l}
\Gamma(t)~\doteq~\ds \Bigg\{ \bfv:[0,T]\mapsto \R^3\,;
~~\hbox{there exists a  positive measure $\mu$, supported on} 
\\[4mm ]\qquad\qquad\qquad \hbox{ the coincidence set 
$\chi(t)$ in (\ref{chit}),
such that for every $s\in [0, T]$ one has}\\[4mm]
\ds  \bfv(s)=-\int_0^s e^{-\beta(t-\sigma)}\left( \int_{[\sigma,t]}  
\Big(\bfn(t,s')\times 
\bigl(\gamma(t,s')-\gamma(t,\sigma)\bigr)\Big)   d\mu(s') \right) \times 
\bigl(\gamma(t,s)-\gamma(t,\sigma)\bigr)\, d\sigma 
\Bigg\}.
\enda
\eeq
Here and in the sequel,  $\bfn(t,s')$ denotes the unit outer normal vector to the 
set $\Omega$ at the boundary point $\gamma(t,s')\in \partial \Omega$.
\v
{\bf Remark 1.} As in  \cite{BPS}, the definition of the cone $\Gamma(t)$ 
in (\ref{RC}) 
is motivated by the following considerations.   
At any point $P=\gamma(t,s')\in \chi(t)$ 
where the stem touches the obstacle, an outward pointing  force acting on the stem 
at $P$ can produce an infinitesimal deformation
described by
$$\gamma^\ve(t,s)~=~\gamma(t,s) + \ve \bfv(s),$$
with
\bel{vvv}\bfv(s)~=~\int_0^s \omega(\sigma)\times 
\bigl(\gamma(t,s)-\gamma(t,\sigma)
\bigr))\, d\sigma\,.\eeq
Here $\omega(\sigma)$ describes the infinitesimal bending of the stem at 
the point $\gamma(t,\sigma)$.   
The elastic energy of the corresponding deformation can be described as
\bel{energy}
{\cal E}(\omega)~=~{1\over 2} \int_0^t e^{\beta(t-\sigma)} |\omega(\sigma)|^2\, d\sigma.\eeq
Notice that the 
weight $e^{\beta(t-\sigma)}$ accounts for the fact that older cells are stiffer, and 
offer more resistance to bending.
It is natural to choose $\omega$ in order to minimize
the total energy ${\cal E}$, subject to a linear constraint of the form
$$\bfn(t,s')\cdot \bfv(s') ~=~c_0$$
for some $c_0>0$.  Necessary conditions for optimality yield the representation
\bel{ooo}\omega(\sigma)~=~\left\{\bega{cl} -\lambda e^{-\beta(t-\sigma)} \bfn(t, s') \times 
\bigl(\gamma(t, s') - \gamma(t,\sigma)\bigr)\qquad&\hbox{for}\quad 0\leq \sigma\leq s',\cr
0\qquad&\hbox{for}\quad s'<\sigma\leq t,\enda\right.\eeq
for some Lagrange multiplier $\lambda>0$.
Inserting (\ref{ooo}) in (\ref{vvv}) and integrating over the set $\chi(t)$ where
the stem touches the obstacle, one formally obtains (\ref{RC}).
\v
%
The equation (\ref{E1}) will be solved on a domain 
of the form 
\bel{Ddef}\D~\doteq~\bigl\{(t,s)\,;~~t\in  [t_0,T]\,,~~s\in [0,T]\bigr\}, \eeq 
with initial and boundary conditions
\bel{iP}
\gamma
(t_0,s)~=~\ov \gamma
(s),\qquad \qquad s\in [0, t_0],\eeq
\bel{bP}
\gamma_{ss}(t,s)~=~0,\qquad\qquad t\in [t_0,T],~~s\in \,]t,T]\,,\eeq
and the  constraint
\bel{Pout}
\gamma
(t,s)~\notin~\Omega\qquad\forall ~t\in [t_0,T],~~s\in [0,t]\,.\eeq
Differentiating w.r.t.~$s$, one obtains an equivalent evolution 
equation for the
unit tangent vector $\bfk$, namely 
\bel{E2}
\bfk_t(t,s)~=~\left(\int_0^{s\wedge t} \Psi\Big(t,\sigma, \gamma(t,\sigma), 
\gamma_s(t,\sigma)\Big)
d\sigma\right)\times  \bfk(t,s)
+\bfh(t,s).
\eeq
Here
$\bfh(t,\cdot)$ 
is any element of the cone
\bel{RC'}\bega{rl}
\Gamma'(t)&\doteq~\ds \Bigg\{ \bfh:[0,t]\mapsto \R^3\,;~~
\hbox{there exists a  positive measure $\mu$ supported on $\chi(t)$
such that}\\[4mm]
&\ds  \bfh(s)~=~ - \int_0^{s}\left(\int_{[\sigma,t]}  e^{-\beta(t-\sigma)}
\bfn(t,s')\times 
\bigl(\gamma
(t,s')-\gamma
(t,\sigma)\bigr)\,d\mu(s')\right) d\sigma\times 
\bfk(t,s)\Bigg\}.
\enda
\eeq
The equation (\ref{E2}) should be solved on the domain $\D$ in (\ref{Ddef}), 
with initial and boundary conditions
\bel{ik}\bfk(t_0, s)~=~\ov \bfk(s)~=~\ov\gamma_s(s),\qquad \qquad s\in [0, t_0]\,,\eeq
\bel{bk}\bfk_s(t,s)~=~0,\qquad\qquad 
t\in [t_0,T],~~s\in [t,T]\,,\eeq
together  with the state constraint  (\ref{Pout}).
Notice that the right hand side of (\ref{E2}) is always perpendicular
to the tangent vector $\bfk(t,s)\doteq \gamma
_s(t,s)$. 
As a consequence, the identities
$$|\bfk(t,s)|~ =~ |\gamma
_s(t,s)|~ =~ 1$$
remain always valid, provided they hold at the initial time $t=t_0$. 
\v
{\bf Definition 1.}  {\it 
We say that a function $\gamma
=\gamma
(t,s)$, defined for 
$(t,s)\in [t_0, T]\times [0,T]$ is a solution 
to the equation (\ref{E1})-(\ref{RC})
with initial and boundary  conditions (\ref{iP})--(\ref{Pout}) if the following holds.
\begi
\item[(i)] The map
$t\mapsto \gamma
(t,\cdot)$ is Lipschitz continuous from $[t_0,T]$ into $H^2([0,T];\,\R^3)$.
\item[(ii)]  For every $t,s$ one has
\bel{inte}\bega{rl}
\gamma
(t,s)&\ds=~\gamma
(t_0,s)+ \int_{t_0}^t \int_0^{s\wedge t} \Psi\Big(\tau,\sigma, \gamma
(\tau,\sigma), 
\gamma
_s(\tau,\sigma)\Big)\times 
\Big(\gamma
(\tau,s)-\gamma
(\tau,\sigma)\Big)
d\sigma\, d\tau\\[4mm]
&\qquad\quad \ds+\int_0^t\bfv(\tau,s)\, d\tau\,,\enda\eeq
where each $\bfv(\tau,\cdot)$ is an element of the cone $\Gamma(\tau)$ defined as in 
(\ref{RC}).
\item[(iii)] The initial conditions hold:
\bel{ic4}
\gamma
(t_0,s)~=~\left\{\bega{cl}\ov \gamma
(s)&\qquad\hbox{if}~~s\in [0, t_0],\\[4mm]
\ov \gamma
(t_0)+(s-t_0)\ov \gamma
'(t_0)&\qquad\hbox{if}~~s\in [t_0,T].
\enda\right.\eeq
\item[(iv)]
The pointwise constraints hold:
\bel{pc4}
 \gamma
(t,s)~\notin~\Omega\qquad\quad\forall t\in [t_0,T],~~s\in [0,t].\eeq
\bel{pc44}
\gamma
(t,s)~=~\gamma(t,t)+(s-t)\gamma_s(t,t)\qquad\quad\forall t\in [t_0,T],~~s\in [t,T].\eeq
 \endi
}

\begin{figure}[ht]
\centerline{\hbox{\includegraphics[width=14cm]{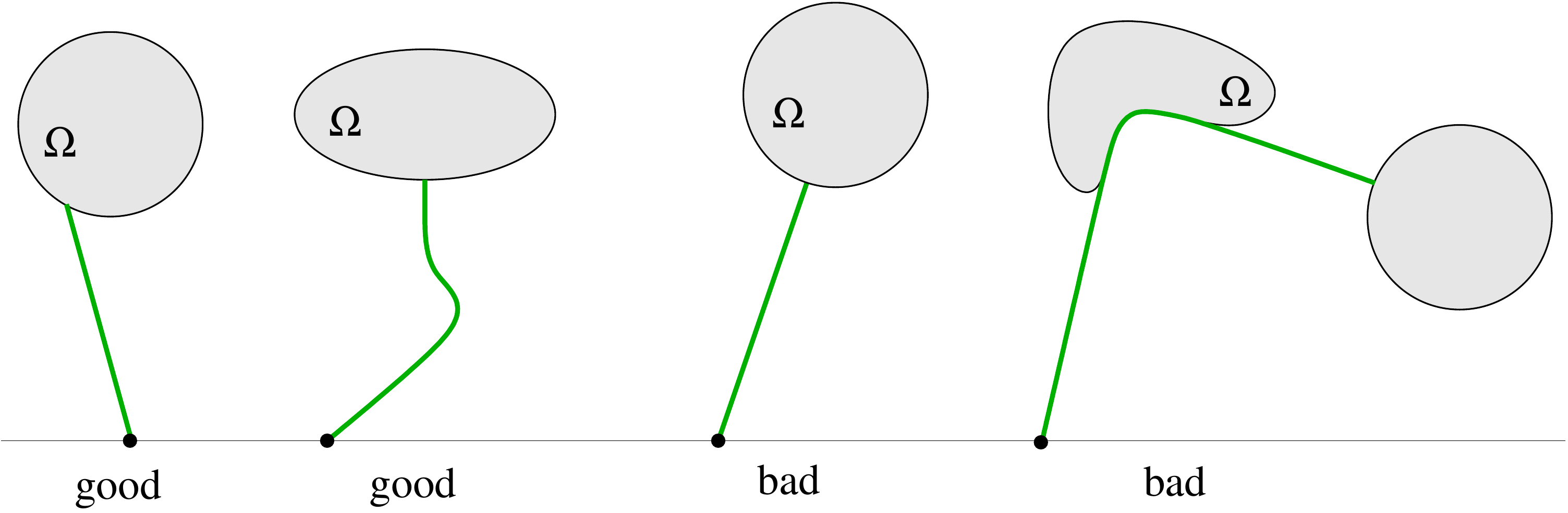}}}
\caption{\small For the two initial configurations on the 
left,  the constrained growth equation 
(\ref{E1}) admits a unique solution. 
On the other hand, the 
two configurations on the right satisfy both (\ref{bad1}) 
and (\ref{bad2}) in
{\bf (B)}. In such cases, the Cauchy problem is ill posed. }
\label{f:sg101}
\end{figure}

Notice that the conditions (\ref{ic4}) and (\ref{pc44}) imply that 
(\ref{iP})-(\ref{bP}) are satisfied.
Given an initial data $\gamma
(t_0,s)=\ov \gamma
(s)$, 
the result in \cite{BPS} provides the existence of a solution 
as long as the following breakdown configuration
is not attained
(see Fig.~\ref{f:sg101}).
\begi
\item[{\bf (B)}] {\it The tip of the stem touches the obstacle perpendicularly, namely
\bel{bad1} \ov \gamma
(t_0)~\in~ \partial \Omega\,,\qquad\qquad 
\ov \gamma
_s(t_0)~=~-\bfn(\ov \gamma
(t_0)).\eeq 
Moreover, 
\bel{bad2} \ov \gamma
_{ss}(s)~= 0\qquad  
\hbox{for all $ s\in ~]0,t[\,$ such that~}\ov \gamma
(s)\notin \partial\Omega\,.\eeq
}
\endi
Here  $\bfn(x)$ denotes the unit outer normal to $\Omega$ at
a boundary point $x\in \partial\Omega$.  

\v
{\bf Theorem 1.} 
{\it Let $\Psi$ in (\ref{E1}) be a $\C^2$ function,
and let
$\Omega\subset\R^3$ be a bounded open set with $\C^2$ boundary. 
At time $t_0$, consider the initial data (\ref{iP}), where
the curve $s\mapsto \ov \gamma
(s)$ is in 
$H^2([0,t_0];\,\R^3)$ and satisfies
\bel{ovP} \ov \gamma
(0)=0\notin \partial\Omega,\qquad\qquad \ov \gamma
(s)\notin\Omega\qquad\forall s\in [0, t_0].\eeq
Moreover, 
assume that the condition {\bf (B)} does NOT hold.

Then there exists $T>t_0$ such that  
the equations (\ref{E1})-(\ref{RC})
with initial and boundary  conditions (\ref{iP})--(\ref{Pout})
admit at least one solution for $t\in [t_0,T]$.

Either (i) the solution is globally defined for all
times $t\geq t_0$, or else (ii) the solution can be extended to
a maximal time interval $[0,T]$, where $\gamma
(T,\cdot)$ satisfies all
conditions in {\bf (B)}.} 
\v
In the present paper we prove that the above solution is unique.
Moreover,
for a.e.~time $t$  the velocity 
$\bfv(t,\cdot)$  determined by the constraint 
reaction
can be computed as follows.
Using the shorter notation 
$\Psi(\sigma)= \Psi\Big(t,\sigma, \gamma
(t,\sigma), 
\gamma
_s(t,\sigma)\Big)$ and $\bfn(t,s)= \bfn(\gamma
(t,s))$ whenever
$\gamma
(t,s)\in 
\partial \Omega$,
consider the minimization
problem
\bel{Min1}
\hbox{minimize:}\quad \E(\omega)~\doteq~{1\over 2}\int_0^te^{\beta(t-\sigma)} |\omega(\sigma)|^2
\, d\sigma,\eeq
subject to the unilateral constraint
\bel{Min2}
\left\langle
\int_0^s\bigl(\Psi(\sigma)+\omega(\sigma)\bigr)\times
\bigl(\gamma
(t,s)-\gamma
(t,\sigma)\bigr)\, d\sigma\,,~\bfn(t,s)\right\rangle
~\geq~0\qquad\qquad\forall s\in \chi(t).
\eeq
If the tip of the stem touches the obstacle, then we  also impose that it
does not penetrate, namely
\bel{Min3}
\left\langle \gamma
_s(t,t)+
\int_0^t\bigl(\Psi(\sigma)+\omega(\sigma)\bigr)\times
\bigl(\gamma
(t,t)-\gamma
(t,\sigma)\bigr)\, d\sigma\,,~\bfn(t,t)\right\rangle
~\geq~0.\eeq
We will show that, at a.e. time $t$, the evolution equation (\ref{E1})
is satisfied with
\bel{E8}
\gamma
_t(t,s)~=~\int_0^s \Big(\Psi(\sigma)+\bar\omega(t,\sigma)\Big)\times 
\Big(\gamma
(t,s)-\gamma
(t,\sigma)\Big)
d\sigma,\eeq
where $\bar\omega(\cdot)$ is the unique minimizer for (\ref{Min1})--(\ref{Min3}).
In other words, for a.e.~time $t$, among all possible choices of 
$\bfv\in\Gamma(t)$, the equation (\ref{E1}) is satisfied precisely with
\bel{Min4}
\bfv(t,s)~=~\int_0^s\bar \omega(t,\sigma)\times
\bigl(\gamma(t,s)-\gamma
(t,\sigma)\bigr)\, d\sigma\,.\eeq
\v
{\bf Theorem 2 (uniqueness).} {\it In the same setting as Theorem~1, the solution
to the evolution equation (\ref{E1})-(\ref{RC}) with initial and boundary 
conditions  (\ref{iP})--(\ref{Pout})  is unique. }
\v
{\bf Theorem 3 (representation of solutions).} {\it
For a.e.~$t\in [0,T]$
the time derivative 
$\gamma_t$ of the solution constructed in Theorem~1 is given by (\ref{E8}),
where $\bar\omega(t,\cdot)$ is the unique minimizer of 
(\ref{Min1}), subject to (\ref{Min2})-(\ref{Min3}).
}
\v
\section{Preliminary lemmas}
\label{s:4}
\setcounter{equation}{0}

In the following,  given a vector $\bfw=(w_1,w_2,w_3)^T$, we shall
denote by $R[\bfw]$  the $3\times 3$ rotation matrix 
\bel{R} R[\bfw]~\doteq~e^A~\doteq~\sum_{k=0}^\infty {A^k\over k!},
\qquad\qquad A~
\doteq~\left(\bega{ccc}
0&-w_3&-w_2\cr w_3&0&-w_1\cr
-w_2&w_1&0\enda\right).\eeq
Notice that, for every $\bar \bfv\in \R^3$, the image
$R[\bfw]
\bar \bfv$ is the value at time $t=1$ of the solution to
$$\dot \bfv(t)~=~\bfw
\times \bfv(t),\qquad\qquad \bfv(0)=\bar \bfv.$$

Next, consider two time-dependent unit vectors:
$\bfk_1(t), \bfk_2(t)$.    
We seek rotation vectors $\bfw(t)$ such that
\bel{rot}
\bfk_2(t)~=~R[\bfw(t)] \,\bfk_1(t)\qquad\quad\forall t\geq 0.\eeq
In particular,  we seek an equation relating
the time derivatives $\bfw_t$ and $\bfk_{i,t}$,
$i=1,2$.
Differentiating (\ref{rot}) w.r.t.~time, one obtains
\bel{rott}
\bfk_{2,t}(t)~=~\left({d\over d \bfw} R[\bfw(t)]\bfw_t\right)
\bfk_1(t)+R[\bfw(t)] \,\bfk_{1,t}(t).\eeq
Assume that
\bel{kder}
\bfk_{i,t}(t)~=~\omega_i(t)\times\bfk_i(t),\qquad\qquad i=1,2,\eeq
for some angular velocities $\omega_1,\omega_2$.
At a time $\tau$ where $\bfw(\tau)=0$, and hence 
$R[\bfw(\tau)] = I$ is the identity matrix, (\ref{rott})
reduces to 
\bel{91}
\omega_2(\tau)\times \bfk_2(\tau)~=~\bfw_t(\tau)\times \bfk_1(\tau) 
+ \omega_1(\tau)\times \bfk_1(\tau).
\eeq
Hence, since $\bfk_2(\tau)= \bfk_1(\tau)$, one has 
$\bfw_t = \omega_2-\omega_1$.
We now study  the more general case where $\bfw(\tau) $ is small 
but nonzero.
\v
{\bf Lemma 1.} {\it   Assume that the  unit vectors 
$\bfk_1(\cdot),\bfk_2(\cdot)$ satisfy
(\ref{kder})
for some continuous angular velocities  $\omega_i(\cdot)$.
Moreover, assume that, at some time $\tau$, one has 
$$\bfk_2(\tau)~=~ R[\bfw(\tau)]\,\bfk_1(\tau),$$
with  $|\bfw(\tau)|<\delta$ sufficiently small.  Then  there exists
$T>\tau$, a constant $C$,  and an absolutely continuous 
map $t\mapsto\bfw(t)$ such that (\ref{rot}) holds 
for all $t\in [\tau,T]$, and moreover
\bel{es1}
\Big|{d\over dt}\bfw(t) - (\omega_2(t)- \omega_1(t)) \Big|~\leq~
C\cdot\Big( |\omega_1(t)|+| \omega_2(t)|\Big)\, |\bfw(t)|.
\eeq
}
\v
{\bf Proof.}
For a fixed $\tau$, choose two additional vectors
$\bfv_1,\bfv_2$  so that $\{ \bfk_2(\tau), \bfv_{1}, \bfv_{2}\}$  
is a (positively oriented) orthonormal basis of $\R^3$. 
Consider the function\footnote{To differentiate the exponential matrix 
$R[\cdot]$, we use the formula ${d\over d\epsilon}
e^{A+\epsilon B}\bigl|_{\epsilon=0}~=~\int_0^1 e^{(1-\xi)
 A} B\,e^{\xi A}\, d\xi
$.}
\bel{F}\bega{rl}
F(c_{1},c_{2})&\ds\doteq ~{d\over dt} \bigg\{
R\Big[\bfw(\tau)+(t-\tau)\bigl(\omega_2(\tau)
-\omega_1(\tau) +c_1\bfv_1 + c_2\bfv_2\bigr)\Big]
\bfk_1(t) -\bfk_2(t) \bigg\}_{t=\tau}\\[4mm]
&=~\ds  \int_0^1 
R\bigl[(1-\xi)
 \bfw(\tau)\bigr]\Big(
\bigl(\omega_2(\tau)
-\omega_1(\tau) +c_1\bfv_1 + c_2\bfv_2\bigr)\times
R\bigl[\xi \bfw(\tau)\bigr] \bfk_1(\tau)\Big) \, d\xi  \\[4mm]
&\qquad\qquad +R[\bfw(\tau)]\,\bigl(\omega_1(\tau)\times\bfk_1(\tau)\bigr)
-\omega_2(\tau)\times \bfk_2(\tau) \\[4mm]
&=~\ds \int_{0}^{1}\bigl(R\bigl[(1-\xi)\bfw(\tau)\bigr](\omega_2(\tau)
-\omega_1(\tau) +c_1\bfv_1 + c_2\bfv_2\bigr)\bigr)d\xi \times \bfk_{2}(\tau)
\\[4mm]
&\qquad\qquad+\bigl(R\bigl[\bfw(\tau)\bigr]\omega_{1}(\tau)-\omega_{2}(\tau)\bigr)\times \bfk_{2}(\tau).
\enda
\eeq
Notice that the vector $F(c_{1},c_{2})$ is always perpendicular to 
$\bfk_{2}(\tau)$.  Hence the vector equation 
\bel{F0}F(c_1,c_2)~=~0\in\R^3\eeq
is equivalent to the  system of two scalar equations
\bel{F12}
F_1(c_1,c_2)~=~\bfv_1\cdot F(c_1,c_2)~=~0,\qquad\qquad
F_2(c_1,c_2)~=~\bfv_2\cdot F(c_1,c_2)~=~0,\eeq
where the dot indicates a scalar product.
The partial derivatives of the map $(c_1,c_2)\mapsto (F_1,F_2)$ are computed by
\bel{derG}
\bega{ccc}\ds
\frac{\partial F_{i}}{\partial c_{j}}~=~
\frac{\partial}{\partial c_{j}}\left(\bfv_i\cdot\int_{0}^{1}\bigl(R\bigl[\xi\bfw(\tau)\bigr](\omega_2(\tau)
-\omega_1(\tau) +c_1\bfv_1 + c_2\bfv_2\bigr)\bigr)\,d\xi \times \bfk_{2}(\tau) \right)
\\[4mm] \ds
=~\bfv_i\cdot 
\int_{0}^{1}\Big(R\bigl[\xi\bfw(\tau)\bigr]\bfv_{j}\Big)d\xi\,\times 
\bfk_2(\tau) ~=~(\bfv_{j}\times \bfk_{2}(\tau))\cdot \bfv_{i}+\mathcal{O}(1)\cdot |\bfw(\tau)|\,.
\enda
\eeq
Hence the Jacobian matrix is
\bel{JF}\left(\frac{\partial F_{i}}{\partial c_{j}}\right)_{i,j=1,2}~=~\left(\bega{cc}
0&1\cr -1 & 0\enda\right)+\mathcal{O}(1)\cdot |\bfw(\tau)|\,.
\eeq
As usual, here the Landau symbol $\O(1)$ denotes a uniformly bounded quantity.
In particular, for $|\bfw(\tau)|$ sufficiently small this Jacobian matrix is invertible. 

We now observe that the right hand side of (\ref{F}) is linear w.r.t.~the vectors
$\omega_1(\tau), \omega_2(\tau)$.    Moreover:
\begi
\item[(i)] When $\bfw(\tau)=0$ we have $\bfk_2(\tau)=\bfk_1(\tau)$ and 
$R[\bfw(\tau)] = I$. In this case, for arbitrary $\omega_1(\tau),\omega_2(\tau)\in\R^3$,
the equation (\ref{F0}) is satisfied by taking $c_1=c_2=0$.
\item[(ii)] When $\omega_1(\tau) = \omega_2(\tau)=0\in\R^3$, for an arbitrary 
$\bfw(\tau)$ the equation 
(\ref{F0}) is  again satisfied by taking $c_1=c_2=0$.\endi

By an application of the implicit function theorem, we obtain the existence of a unique vector, say
\bel{oms}\omega^\sharp(\tau)~=~c_1\bfv_1 + c_2\bfv_2~=~ \Phi\bigl(\bfk_2(t), \bfw(t), 
\omega_1(t),\omega_2(t)\bigr),\eeq
satisfying 
\bel{Eq1}\omega^\sharp(\tau)~\in~\bfk_2(\tau)^\perp,\eeq
\bel{Eq2} 
\bega{l}\ds
\ds \int_{0}^{1}\bigl(R\bigl[(1-\xi)\bfw(\tau)\bigr](\omega_2(\tau)
-\omega_1(\tau) +\omega^{\sharp}(\tau)\bigr)\bigr)d\xi \times \bfk_{2}(\tau)
\\[4mm]
\qquad\qquad+\bigl(R\bigl[\bfw(\tau)\bigr]\omega_{1}(\tau)-\omega_{2}(\tau)\bigr)\times \bfk_{2}(\tau)~=~0.
\enda
\eeq
The above identities (i)-(ii) imply
\bel{Hom}
|\omega^\sharp(\tau)|~=~
 \O(1)\cdot \Big( |\omega_1(\tau)|+| \omega_2(\tau)|\Big)\, |\bfw(\tau)|.\eeq 
By the continuity of the angular velocities $\omega_1,\omega_2$, the 
above construction can be repeated for every $t\in [\tau,T]$, as long as
the rotation vector $\bfw(t)$ remains sufficiently small. 
This yields an evolution equation for $\bfw$, of the form
\bel{wt}{d\over dt}\bfw(t)~=~\omega_2(t)-\omega_1(t) + \Phi\bigl(\bfk_2(t), \bfw(t), 
\omega_1(t),\omega_2(t)\bigr),\eeq
where $\Phi$ is the function implicitly defined in (\ref{oms}), providing the unique solution to 
(\ref{Eq1})-(\ref{Eq2}).
By the regularity of $\Phi$, 
given the functions $\omega_1(\cdot)$, $\omega_2(\cdot)$,  $\bfk_2(\cdot)$ 
and the initial condition
$\bfw(\tau)$,  the evolution equation (\ref{wt}) has a unique local solution, defined as long as the vector $\bfw$ remains small enough.
This completes the proof of the lemma.
\endproof


\v
Toward a proof of Theorem 2 we need an integral version of
Lemma~1.   As before, we  consider two curves, growing in time:
$\gamma_i(t,s)$, $s\in [0,t]$.   We denote by $\bfk_i(t,s)=\gamma_{i,s}(t,s)$
the unit tangent vectors.
\v
{\bf Lemma 2.}  {\it
Assume  that, for $i=1,2$, 
\bel{12}
\bfk_{i,t}(t,s)~=~\left(\int_0^s\omega_i(t,\sigma)\, d\sigma
\right)\times\bfk_i(t,s),\qquad\qquad s\in [0,t].\eeq
Moreover, assume that at time $\tau$ one has
\bel{gtau}\bfk_2(\tau,s)~ =~ R\left[\int_0^s\bfw(\tau,\sigma)d\sigma\right]
\bfk_1(\tau,s),\eeq
with $\|\bfw(\tau,\cdot)\|_{\L^{2}([0,\tau])}\leq\delta$ sufficiently small.
Then there exists $T>\tau$ such that, for all $t\in [\tau,T]$ one has
the representation 
\bel{gt}\bfk_2(t,s)~ =~ R\left[\int_0^s\bfw(t,\sigma)d\sigma\right]
\bfk_1(t,s),\qquad\qquad s\in [0, t].\eeq
Here the rotation vectors $\bfw(t,\cdot)$ can
be chosen  so that 
\bel{dom12}\ds\Big|\int_{0}^{s}\bfw_t(t,\sigma)-\omega_2(t,\sigma)+\omega_1(t,\sigma)d\sigma
\Big|~=~\O(1)\cdot 
 \Big( \Big|\int_{0}^{s}|\omega_1(t,\sigma)d\sigma\Big|+\Big| \int_{0}^{s}\omega_2(t,\sigma)d\sigma \Big|\Big)\, 
 \|\bfw(t,\cdot)\|_{\L^1([0,t])}\,\eeq
 }
 for $s\in[0,t]$.
 \v
 {\bf Proof.}   We repeat the construction of Lemma~1.
For every $s\in [0,\tau]$ we have
 $$|\bfk_1(\tau,s)|~=~|\bfk_2(\tau,s)|~=~1,
\qquad\quad \bigl|\bfk_2(\tau,s)-\bfk_1(\tau,s)\bigr|~=~\O(1)\cdot \|\bfw(\tau,\cdot)
\|_{\L^1([0,\tau])}\,.$$
For each $s\in [0,\tau]$, choose unit vectors $\bfv_1(s), \bfv_2(s)$ so that 
$\{\bfk_2(\tau,s), \bfv_1(s),\bfv_2(s)\}$ is a (positively oriented)
orthonormal basis of $\R^3$. Notice that $s\mapsto \bfv_{1}(s), \bfv_{2}(s)$ are in $H^{1}([0,\tau])$.

Given two scalar functions $c_1(s)$, $c_2(s)$,
for each $s\in [0,\tau]$ define
{\small
\bel{FF2}\bega{l}\ds
F(s,c_{1}(s),c_{2}(s))\\[4mm]
\doteq \ds~{d\over dt} \bigg\{
R\left[\int_0^s\bfw(\tau,\sigma)+(t-\tau)\bigl(\omega_2(\tau,\sigma)
-\omega_1(\tau,\sigma) +c_1(\sigma)\bfv_1(\sigma) + c_2(\sigma)\bfv_2(\sigma)
\bigr)\, d\sigma\right]
\bfk_1(t,s) -\bfk_2(t,s) \bigg\}_{t=\tau}\\[4mm]
=~\ds  \int_0^1 
R\left[(1-\xi)\int_0^s
 \bfw(\tau,\sigma)\, d\sigma\right]\left(\int_0^s
\bigl(\omega_2(\tau,\sigma)
-\omega_1(\tau,\sigma) +c_1(\sigma)\bfv_1(\sigma) + c_2(\sigma)\bfv_2(\sigma)\bigr)
d\sigma\right)\\[4mm]
\ds \qquad\qquad \times
R\left[\xi \int_0^s\bfw(\tau,\sigma)\, d\sigma\right] \bfk_1(\tau,s)
 \, d\xi  \\[4mm]
\ds\qquad\qquad +R\left[\int_0^s\bfw(\tau,\sigma)\, d\sigma\right]\,\left(\int_0^s\omega_1(\tau,\sigma)\, d\sigma\right)\times\bfk_1(\tau,s)
-\left(\int_0^s\omega_2(\tau,\sigma)\, d\sigma\right)\times \bfk_2(\tau,s)
\\[4mm]
=~\ds   \int_0^1 
R\left[\xi\int_0^s
 \bfw(\tau,\sigma)\, d\sigma\right]\left(\int_0^s
\bigl(\omega_2(\tau,\sigma)
-\omega_1(\tau,\sigma) +c_1(\sigma)\bfv_1(\sigma) + c_2(\sigma)\bfv_2(\sigma)\bigr)
d\sigma\right) d\xi \times \bfk_{2}(\tau,s)
\\[4mm]\ds
\qquad\qquad+\left(R\left[\int_0^s
 \bfw(\tau,\sigma)\, d\sigma\right]
  \int_0^s \omega_{1}(\tau,\sigma) \, d\sigma\right) \times \bfk_{2}(\tau,s)- \int_0^s \omega_2(\tau,\sigma) \, d\sigma \times \bfk_{2}(\tau,s).
\enda
\eeq
}
Notice that the vector $F(s,c_{1}(s),c_{2}(s))$ is always perpendicular to 
$\bfk_{2}(\tau,s)$.  Hence the vector equation 
\bel{F3}F(s, c_1(s),c_2(s))~=~0\in\R^3\eeq
is equivalent to the  system of two scalar equations
\bel{F4}\left\{
\bega{rl}F_1(s,c_1(s),c_2(s))&\doteq
~\bfv_1(s)\cdot F(s,c_1(s),c_2(s))~=~0,\\[3mm]
F_2(s,c_1(s),c_2(s))&\doteq~\bfv_2(s)\cdot F(s,c_1(s),c_2(s))~=~0.
\enda\right.\eeq
These should hold for all $s\in [0,\tau]$.

For $s=0$ the equations (\ref{F4}) are trivially satisfied.
Hence it suffices to solve the equations for the derivatives:
\bel{F5} {d\over ds}F_i(s,c_1(s),c_2(s))~=~0,\qquad\qquad i=1,2,\qquad s\in [0,\tau].\eeq

Notice also that $$\bfv_{1}\cdot(\bfa\times \bfk_{2})=\bfa\cdot(\bfk_{2}\times \bfv_{1})=\bfv_{2}\cdot \bfa, $$
$$\bfv_{2}\cdot(\bfa\times \bfk_{2})=\bfa\cdot(\bfk_{2}\times \bfv_{2})=-\bfv_{1}\cdot \bfa,$$

and that, in view of \eqref{F3}, 
\bel{F_eps}
\ds{d\over ds} F_{i}(s,c_{1}(s),c_{2}(s))=\Big( {d\over ds} F(s,c_{1}(s),c_{2}(s))\Big)\cdot \bfv_{i}(s).
\eeq
Taking these observations into account, we then compute:
\bel{F6} \bega{l}\ds{d\over ds}F_1(s,c_1(s),c_2(s))
\\[4mm]
\ds=~\int_{0}^{1} \int_{0}^{1} \Big( \alpha\, \xi \,R\Big[\int_{0}^{s}\bfw(\tau,\sigma)d\sigma \Big]\cdot \bfw(\tau,s)\Big)
\\[4mm]
\qquad \qquad\ds \times ~\Big(\int_{0}^{s}(\omega_{2}(\tau,\sigma)-\omega_{1}(\tau,\sigma)+c_{1}(\sigma)\bfv_{1}(\sigma)+c_{2}(\sigma)\bfv_{2}(\sigma))d\sigma\Big)d\alpha\, d\xi \cdot \bfv_{2}(s)
\\[4mm]
\qquad\ds+~\int_{0}^{1} R\Big[ \xi \int_{0}^{s}\bfw(\tau,\sigma)d\sigma\Big](\omega_{2}(\tau,s)-\omega_{1}(\tau,s) +c_{1}(s)\bfv_{1}(s)+c_{2}(s)\bfv_{2}(s))\,d\xi \cdot \bfv_{2}(s)
\\[4mm]
\qquad\ds+~ \int_{0}^{1}\Big(R\Big[\xi\, \int_{0}^{s}\bfw(\tau,\sigma)d\sigma \Big]\cdot \bfw(\tau,s)\,\times \Big(\int_{0}^{s}\omega_{1}(\tau,\sigma)\, d\sigma \Big)\Big)d\xi \cdot \bfv_{2}(s)
\\[4mm]
\qquad\ds+~\Big( R\Big[ \int_{0}^{s}\bfw(\tau,\sigma)d\sigma\Big]\cdot \omega_{1}(\tau,s) \Big)\cdot \bfv_{2}(s) - \omega_{2}(\tau,s)\cdot \bfv_{2}(s)~=~0.
\enda \eeq

A similar relation holds true for ${d\over ds}F_{2}(s,c_{1}(s),c_{2}(s)).$ Notice that both relations together lead to a system of equations
\bel{C1}
\bega{ll}\ds c_1(s)~=~-\int_{0}^{1} \int_{0}^{1} \Big( \alpha\, \xi \,R\Big[\int_{0}^{s}\bfw(\tau,\sigma)d\sigma \Big]\cdot \bfw(\tau,s)\Big)
\\[4mm]
\qquad \qquad\qquad \ds \times ~\Big(\int_{0}^{s}(\omega_{2}(\tau,\sigma)-\omega_{1}(\tau,\sigma)+c_{1}(\sigma)\bfv_{1}(\sigma)+c_{2}(\sigma)\bfv_{2}(\sigma))d\sigma\Big)d\alpha\, d\xi \cdot \bfv_{1}(s)
\\[4mm]
\qquad \qquad\ds + ~ \Big(\mathcal{O}(1)\cdot\Big|\int_{0}^{s}\bfw(\tau,\sigma)d\sigma\Big|(\omega_{2}(\tau,s)-\omega_{1}(\tau,s) +c_{1}(s)\bfv_{1}(s)+c_{2}(s)\bfv_{2}(s))\, \Big) \cdot \bfv_{1}(s)
\\[4mm]
\qquad \qquad\ds + \Psi_{1}(s,\omega_{1},\omega_{2},\bfw, \bfv_{1}) 
~\doteq~ P_{1}(\bfc)
\enda\eeq
\bel{C2}
\bega{ll}\ds c_2(s)~=~ - \int_{0}^{1} \int_{0}^{1} \Big( \alpha\, \xi \,R\Big[\int_{0}^{s}\bfw(\tau,\sigma)d\sigma \Big]\cdot \bfw(\tau,s)\Big)
\\[4mm]
\qquad \qquad\qquad \ds \times ~\Big(\int_{0}^{s}(\omega_{2}(\tau,\sigma)-\omega_{1}(\tau,\sigma)+c_{1}(\sigma)\bfv_{1}(\sigma)+c_{2}(\sigma)\bfv_{2}(\sigma))d\sigma\Big)d\alpha\, d\xi \cdot \bfv_{2}(s)
\\[4mm]
\qquad \qquad\ds + ~ \Big(\mathcal{O}(1)\cdot\Big|\int_{0}^{s}\bfw(\tau,\sigma)d\sigma\Big|(\omega_{2}(\tau,s)-\omega_{1}(\tau,s) +c_{1}(s)\bfv_{1}(s)+c_{2}(s)\bfv_{2}(s))\, \Big) \cdot \bfv_{2}(s)
\\[4mm]
\qquad \qquad\ds + \Psi_{2}(s,\omega_{1},\omega_{2},\bfw, \bfv_{2}) 
~\doteq ~P_{2}(\bfc)
\enda\eeq
where $\Psi_{i}$ for $i=1,2,$  
are  smooth functions which do not depend on $c_{i}(s)$. 
Now denote with $\bfc = (c_{1},c_{2})$ and consider the operator 
$\mathcal{P}[\bfc]\doteq (P_{1}(\bfc),P_{2}(\bfc))$ such that 
$\bfc=\mathcal{P}[\textbf{\bfc}]$. We now aim to show that the 
system \eqref{C1}, \eqref{C2} admits a unique solution proving that 
$\mathcal{P}[\cdot]$ is a contraction on $\L^{2}[0,\tau]$ 
for $\delta$ small enough. Indeed, 
for any $\bfc,\tilde{\bfc}\in \L^{2}[0,\tau]$, one has
\bel{FP}\bega{lll}\ds \bigl\|\mathcal{P}[\bfc]-\mathcal{P}[\tilde{\bfc}]
\bigr\|_{\L^{2}[0,\tau]}~\leq ~ \int_{0}^{\tau}  \Big| 
(I+\mathcal{O}(1)\cdot\delta)\, \bfw(\tau,s)\Big|^{2}
\\[4mm]
\qquad \qquad\qquad\ds \times~
\Big|\int_{0}^{s} (c_{1}(\sigma)-\tilde{c}_{1}(\sigma))\bfv_{1}(\sigma)+(c_{2}(\sigma)-\tilde{c}_{2}(\sigma))\bfv_{2}(\sigma)d\sigma\Big|^{2} 
\\[4mm]
\qquad \qquad\ds +~K_{2}\delta \int_{0}^{\tau} \Big| (c_{1}(\sigma)-\tilde{c}_{1}(\sigma))\bfv_{1}(\sigma)+(c_{2}(\sigma)-\tilde{c}_{2}(\sigma))\bfv_{2}(\sigma)\Big|^{2}d\sigma
\\[4mm]
\qquad \ds \leq ~ (1+K_{1}\delta)^{2} \int_{0}^{\tau}|\bfw(\tau,s)|^{2}ds 
\\[4mm]
\qquad \qquad\qquad \ds \times~
\int_{0}^{\tau} |(c_{1}(\sigma)-\tilde{c}_{1}(\sigma))\bfv_{1}(\sigma)+(c_{2}(\sigma)-\tilde{c}_{2}(\sigma))\bfv_{2}(\sigma)|^{2}d\sigma  
\\[4mm]
\qquad \qquad\ds +~K_{2}\delta \int_{0}^{\tau} \Big| (c_{1}(\sigma)-\tilde{c}_{1}(\sigma))\bfv_{1}(\sigma)+(c_{2}(\sigma)-\tilde{c}_{2}(\sigma))\bfv_{2}(\sigma)\Big|^{2}d\sigma
\\[4mm]
\qquad \ds \leq ~ \delta\Big((1+K_{1}\delta)^{2}+K_{2}\Big)\|\bfc-\tilde{\bfc}\|_{\L^{2}[0,\tau]},
\enda\eeq
for some constants $K_1, K_2$ independent of $\bfc,\tilde{\bfc}$. When $\delta>0$ is sufficiently small, \eqref{FP} shows that $\mathcal{P}[\cdot]$ is a strict contraction. As a consequence, the system of equations \eqref{C1}-\eqref{C2} admits a unique solution, which we denote by  
$\bigl(\bar{c}_{1}(\cdot),\,\bar{c}_{2}(\cdot)\bigr)$. 
Then $\bar{c}_{1}$, $\bar{c}_2$ will also satisfy the 
relations (\ref{F4}). Moreover 
\bel{int_order}
\ds \int_{0}^{s}(\bar{c}_{1}(\sigma)\bfv_{1}(\sigma)+\bar{c}_{2}(\sigma)\bfv_{2}(\sigma))d\sigma~=~\Phi\bigl(\bfk_2(\tau,s), \textbf{W}(\tau,s) , \Omega_1(\tau,s),\Omega_2(\tau,s)\bigr),
\eeq
where 
\bel{int}\ds
\textbf{W}(\tau,s)~=~\int_{0}^{s}\bfw(\tau,\sigma)d\sigma,\qquad \Omega_{1}(\tau,s)=\int_{0}^{s}\omega_{1}(\tau,\sigma)d\sigma,\qquad \Omega_{2}(\tau,s)=\int_{0}^{s}\omega_{2}(\tau,\sigma)d\sigma,
\eeq
and $\Phi$ is a smooth function satisfying
\bel{int_order_cont}\bega{lll}\ds
\Big|\Phi\bigl(\bfk_2(\tau,s), \textbf{W}(\tau,s) , \Omega_1(\tau,s),\Omega_2(\tau,s)\bigr)\Big|~=~\O(1)\cdot\bigl(\|\omega_{1}(\tau,\cdot)\|_{\L^{1}}+
\|\omega_{2}(\tau,\cdot)\|_{\L^{1}}\bigr)\|\bfw(\tau,\cdot)\|_{\L^{1}}
\enda
\eeq
Again, by the continuity of the integrated angular velocities $\Omega_1,\Omega_2$, the 
above construction can be repeated for every $t\in [\tau,T]$, as long as
the rotation vector $\|\bfw(\tau,\cdot)\|_{\L^{2}}$ remains sufficiently small. 
This yields an evolution equation for $\textbf{W}$, of the form
\bel{wt_int}\ds{d\over dt}\textbf{W}(t,s)~=~\Omega_2(t,s)-\Omega_1(t,s) + \Phi\bigl(\bfk_2(t,s), \textbf{W}(t,s), 
\Omega_1(t,s),\Omega_2(t,s)\bigr),\qquad s\in[0,t].\eeq
By the regularity of $\Phi(\cdot,\cdot,\cdot,\cdot)$, 
given the functions $\Omega_1$, $\Omega_2$,  $\bfk_2$, 
and the initial condition $\textbf{W}(\tau,s)$,  the evolution equation (\ref{wt_int}) has a unique local solution for every $s\in[0,t]$, defined as long as the vector $\textbf{W}(t,s)$ remains small enough.
This completes the proof of the lemma. \endproof
\v
\section{Uniqueness of solutions}
\label{s:5}
\setcounter{equation}{0}
Consider two solutions $\gamma_1,\gamma_2$ of (\ref{E1})-(\ref{RC}), and call
$\bfk_i(t,s)=\gamma_{i,s}(t,s)$ the corresponding tangent vectors.
For each $t$, we shall 
construct a rotation vector $\bfw(t,\cdot)$ such that
\bel{kw1}
\bfk_2(t,s)~=~R\left[\int_0^s\bfw(t,\sigma)\, d\sigma\right]\bfk_1(t,s),
\qquad\qquad s\in [0,t]\,.
\eeq
To measure the size of this vector $\bfw$, for any $t>0$ on the space
$\L^2([0,t])$ we shall use the equivalent inner product and norm
\bel{wn}
\langle \bfv,\bfw\rangle~\doteq~\int_0^t e^{-\beta s}
\langle \bfv(s),\,\bfw(s)\rangle\, ds\,,\qquad\qquad \|\bfv\|
~\doteq~\langle \bfv,\bfv \rangle^{1/2}.\eeq
Using this equivalent norm, we shall prove
the key inequality
\bel{w3}
\la \bfw_t(t,\cdot),\, \bfw(t,\cdot)
\ra~\leq~C\cdot \|\bfw(t,\cdot)\|^2,\eeq
for a suitable constant $C$.
In turn, this yields the estimate
\bel{wtn}{d\over dt}\|\bfw(t,\cdot)\|^2~\leq~2C\cdot
\|\bfw(t,\cdot)\|^2\,.\eeq
In particular, if $\bfw(t_0,\cdot)\equiv 0$, 
this will imply $\bfw(t,\cdot)\equiv 0$
for all $t\geq t_0$, proving uniqueness.

Toward a proof of (\ref{w3}) we use the representation
\bel{rep}
\bfk_{i,t}(t,s)~=~\left(\int_0^s \omega_i(t,\sigma)\, d\sigma\right)
\times \bfk_i(t,s),\qquad\qquad s\in [0,t],~~i=1,2,\eeq
where the angular velocities $\omega_i$ satisfy
\bel{omi}\bega{rl}
\omega_i(t,s)&\ds=~\Psi\bigl(t,s,\gamma_i(t,s), 
\gamma_{i,s}(t,s)\bigr)\\[4mm]
&\qquad\ds  - \int_0^{s}\left(\int_{[\sigma,t]}  e^{-\beta(t-\sigma)}
\bfn(t,s')\times 
\bigl(\gamma_i(t,s')-\gamma_i(t,\sigma)\bigr)\,d\mu_i(s')\right) d\sigma\,,\enda
\eeq
where $\mu_i$ is a positive measure, supported on the contact set 
$\{s\in [0, t]\,;~~\gamma_i(t,s)\in \partial\Omega\}$.

Thanks to Lemma 2, since we know that $\omega_1,\omega_2$ 
are uniformly bounded, we have
\bel{osw}
\bigl\langle \bfw_t(t,\cdot)\,,~\bfw(t,\cdot)\bigr\rangle
~=~\bigl\langle \omega_2(t,\cdot)-\omega_1(t,\cdot)\,,~\bfw(t,\cdot)\bigr\rangle+\O(1)\cdot \|\bfw(t,\cdot)\|^2.\eeq
To estimate the first term on the right hand side of (\ref{osw}), we write
\bel{os2}\bega{l}
\bigl\langle \omega_2(t,\cdot)-\omega_1(t,\cdot)\,,
~\bfw(t,\cdot)\bigr\rangle
\\[4mm]
\ds \quad =~\la \Psi\bigl(t,\cdot,\gamma_2(t,\cdot), 
\gamma_{2,s}(t,\cdot)\bigr)-\Psi\bigl(t,\cdot,\gamma_1(t,\cdot), 
\gamma_{1,s}(t,\cdot)\bigr)\,,~\bfw(t,\cdot)\ra\\[4mm]
\ds \qquad+\left\langle \int_0^\cdot\left(\int_{[\sigma,t]}  e^{-\beta(t-\sigma)}
\bfn(t,s')\times 
\bigl(\gamma_1(t,s')-\gamma_1(t,\sigma)\bigr)\,d\mu_1(s')\right) d\sigma\,,
~\bfw(t,\cdot)\right\rangle\\[4mm]
\ds \qquad-\left\langle \int_0^\cdot\left(\int_{[\sigma,t]}  e^{-\beta(t-\sigma)}
\bfn(t,s')\times 
\bigl(\gamma_2(t,s')-\gamma_2(t,\sigma)\bigr)\,d\mu_2(s')\right) d\sigma\,,
~\bfw(t,\cdot)\right\rangle
\\[4mm]
\quad \doteq~J_0+J_1+J_2\,.
\enda\eeq
The regularity properties of $\Psi$ immediately imply
\bel{J0}\bega{rl}|J_0|&\leq~\ds\int_0^t
\Big|\Psi\bigl(t,s,\gamma_2(t,s), 
\gamma_{2,s}(t,s)\bigr)-\Psi\bigl(t,s,\gamma_1(t,s), 
\gamma_{1,s}(t,s)\bigr)\Big|\,|\bfw(t,s)|\, ds
\\[4mm]
&\leq C_0\cdot \|\bfw(t,\cdot)\|^2,\enda\eeq
for some constant $C_0$.

It remains to estimate  the last two terms in (\ref{omi}).
To fix the ideas, consider a point $s'\in \chi_1(t)$, so that 
$\gamma_1(t, s')\in \partial\Omega$.
This point will contribute to the angular velocity $\omega_1$
through a term of the form
\bel{om4}
\left\{
\bega{cl}
e^{-\beta(t-\sigma)} \Big((\gamma_1(t,s')-\gamma_1(t,\sigma))
\times\bfn_1(t,s')\Big)&\qquad\qquad
\hbox{if}~~\sigma\in [0,s'],\\[3mm]
0&\qquad\qquad\hbox{if}~~\sigma>s'.\enda\right.\eeq

By assumption, 
\bel{g2}\bega{l}\ds
\gamma_2(t,s')-\gamma_1(t,s')~
=~\int_0^{s'} \bigl(\bfk_2(t,s)-\bfk_1(t,s)
\bigr)\, ds\\[4mm]
\qquad \ds=~\int_0^{s'} \left(R\left[\int_0^s \bfw(t,\sigma)\, d\sigma
\right] - I\right) \bfk_1(t,s)\, ds \\[4mm]
\qquad\ds =~\int_0^{s'} \left(\int_0^s \bfw(t,\sigma)\, d\sigma
\right)\times \bfk_1(t,s)\, ds +\O(1)\cdot \|\bfw(t,\cdot)\|^2\\[4mm]
\qquad\ds =~\int_0^{s'} \bfw(t,s)\bigl(\gamma_1(t,s') - \gamma_1(t,s)
\bigr)\, ds +\O(1)\cdot \|\bfw(t,\cdot)\|^2\\[4mm]
\qquad =~\O(1)\cdot \|\bfw(t,\cdot)\|^2.
\enda
\eeq
Using (\ref{g2}) and the properties of the triple product, we now compute
\bel{do2}\bega{rl}&\ds
\int_0^{s'}e^{-\beta s}\,
e^{-\beta(t-s)} \la\bigl(\gamma_1(t,s')-\gamma_1(t,s)\bigr)
\times\bfn_1(t,s')\,,~\bfw(t,s)\ra\, ds
\\[4mm]&\ds\qquad
=~e^{-\beta t} \bfn_1(t,s')\cdot\int_0^{s'} \bfw(t,s)\times
\bigl(\gamma_1(t,s')-\gamma_1(t,s)\bigr)\, ds\\[4mm]
&\qquad 
\ds =~e^{-\beta t} \bfn_1(t,s')
\cdot \bigl(\gamma_2(t,s')-\gamma_1(t,s')\bigr)+\O(1)\cdot \|\bfw(t,\cdot)\|^2\\[4mm]
&\qquad\leq ~
\O(1)\cdot \bigl|\gamma_2(t, s') - \gamma_1(t, s')\bigr|^2 +
\O(1)\cdot \|\bfw(t,\cdot)\|^2\\[4mm]
&\qquad =~\O(1)\cdot \|\bfw(t,\cdot)\|^2.
\enda
\eeq
Recalling that the total mass of the measure $\mu_1$ 
is uniformly bounded, the second term on the right hand side of (\ref{os2})
can thus be estimated by
\bel{J1}
J_1~\leq~\O(1)\cdot\,\|\bfw(t,\cdot)\|^2\cdot \int_{[0,t]} \mu_1(s')\, ds'
~\leq~ C_1\,\|\bfw(t,\cdot)\|^2,\eeq
for some constant $C_1$.
Similarly,
\bel{J2}
J_2~\leq~\O(1)\cdot\,\|\bfw(t,\cdot)\|^2\cdot \int_{[0,t]} \mu_2(s')\, ds'
~\leq~C_2\,\|\bfw(t,\cdot)\|^2.\eeq
By (\ref{os2}) together (\ref{J0}), (\ref{J1}), and (\ref{J2}), in view of 
(\ref{osw})  we achieve a proof of (\ref{w3}). 
  By Gronwall's lemma, this proves the 
  uniqueness of solutions to (\ref{E1}) and (\ref{bP})--(\ref{Pout}),
  and continuous dependence of solutions on the initial data (\ref{iP}).
\endproof

\v
\section{Proof of the representation formula}
\label{s:6}
\setcounter{equation}{0}
In this section we give a proof of Theorem 3, showing that 
any  solution to (\ref{E1})-(\ref{RC}) has the form 
(\ref{E8}).

For any time $t$, consider the contact set $\chi(t)$ of points $s\in [0,t]$
where the stem touches the obstacle.   Observe that the map
$t\mapsto \chi(t)$ is an upper semicontinuous multifunction with compact values.

{\bf Lemma 3.}  {\it There exists a set of times $\N$ of measure zero 
such that, for each $t\in [t_0,T]\setminus\N$  the partial  derivative $\gamma_t(t,s)$
exists for all $s\in [0,T]$. Moreover, the map $s\mapsto \gamma_t(\tau,s)$ is 
Lipschitz continuous.} 
\v
{\bf Proof.} 
We  use the representation
$$\gamma(t,s)~=~\int_0^s \bfk(t,\sigma)\, d\sigma,$$
\bel{gt1}\gamma(t+\ve,s) - \gamma(t,s)~=~\int_t^{t+\ve}
\int_0^s \bfk_t(\tau,\sigma)\, d\sigma \, d\tau\,.\eeq
By the regularity of the solution $\gamma$, proved in Theorem 1 of \cite{BPS},
the partial derivative $\bfk_t$ is well defined for a.e.~$(\tau,\sigma)\in [t_0,T]\times[0,T]$. Moreover,  it satisfies a uniform bound $|\bfk_{t}(\tau,\sigma)|\leq C$.

Therefore, 
there exists a set of times $\N\subset [t_0,T]$ of measure zero such that, for 
$t\notin \N$,
the partial derivative $\bfk_t(t,\sigma)$ exists for a.e.~$\sigma\in [0,T]$.

Using (\ref{gt1}) and  the Lebesgue dominated convergence theorem, 
for every $t\notin\N$ we obtain
$$\gamma_t(t,s)~=~\lim_{\ve\to 0}  \frac{\gamma(t+\ve,s) - \gamma(t,s)}{\ve}~=~
\lim_{\ve\to 0} \int_0^s {\bfk(t+\ve,\sigma)- \bfk(t,\sigma)\over\ve}\, d\sigma~=~
\int_0^s\bfk_t(t,\sigma)\, d\sigma.$$
This achieves the proof.
\endproof
\v
{\bf Corollary 1.} {\it Consider any time $\tau\in [t_0,T]\setminus \N$.
Then, calling $\bfn(\tau,s)$ the unit outer normal to the obstacle
at the boundary point $\gamma(\tau,s)\in \partial\Omega$, one has
 \bel{Pt1}
\bigl\langle \gamma_t(\tau,s),\,\bfn(\tau,s)\bigr\rangle~=~0
\qquad\qquad \forall s\in \chi(\tau)\setminus\{\tau\}.\eeq
In addition, if the tip of the stem touches the obstacle,
i.e.~if $\tau\in \chi(\tau)$, then 
\bel{Pt2}
\bigl\langle \gamma_t(\tau,\tau)+ \bfk(\tau,\tau),\,
\bfn(\tau,\tau)\bigr\rangle~=~0.\eeq
}

{\bf Proof.}  Denote by
$$\Phi(x)~\doteq~\left\{\bega{cl}
 \hbox{dist}(x,\Omega)\qquad &\hbox{if}~~x\notin\Omega,\cr
-\hbox{dist}(x,\partial \Omega)\qquad &\hbox{if}~~x\in\Omega,\enda\right.$$
the signed distance of a point $x$ to the boundary of $\Omega$.
Since $\Omega$ has a $\C^2$ boundary, the  function $\Phi$ is twice continuously differentiable in a neighborhood of $\partial \Omega$.

If (\ref{Pt1}) fails for some  $s\in \chi(\tau)\setminus\{\tau\}$,  then
$$\Phi(\gamma(\tau,s))~=~0,\qquad\qquad {d\over dt}\Phi(\gamma(t,s))\bigg|_{t=\tau}~=~
\bigl\langle \gamma_t(\tau,s),\,\bfn(\tau,s)\bigr\rangle~\not= ~0.$$
This yields a contradiction, because $\Phi(\gamma(t,s))~\geq ~0$ for all $t$.

Similarly, if  $\tau\in \chi(\tau)$ but (\ref{Pt2}) fails, then 
$$\Phi(\gamma(\tau,\tau))~=~0,\qquad\qquad {d\over dt}\Phi(\gamma(t,t))\bigg|_{t=\tau}~=~
\bigl\langle \gamma_t(\tau,\tau)+ \bfk_s(\tau,\tau),\,
\bfn(\tau,\tau)\bigr\rangle~\not=~0.$$
This yields a contradiction, because $\Phi(\gamma(t,t))~\geq ~0$ for all $t$.
\endproof
\v
{\bf Proof of Theorem 3.} 

We will show that the representation formula (\ref{E8})
is valid at every 
time $\tau\in [t_0, T]\setminus\N$,  where the conclusions
(\ref{Pt1})-(\ref{Pt2}) of Corollary~1 hold. Notice that, since condition \textbf{(B)} does NOT hold, the set of $\omega$ satisfying the constraints \eqref{Min2}, \eqref{Min3} is non empty.
\v
{\bf 1.} Fix a time  $\tau\in [t_0, T]\setminus\N$ and let 
 $\bfv\in \Gamma(\tau)$ be a velocity field for which the
bilateral  constraints are satisfied:
\bel{Min7}
\left\langle
\int_0^s\Psi(\sigma)\times
\bigl(\gamma
(\tau,s)-\gamma
(\tau,\sigma)\bigr)\, d\sigma + \bfv(s)\,,~\bfn(\tau,s)\right\rangle
~=~0\qquad\qquad\forall s\in \chi(\tau)\setminus \{\tau\},
\eeq
 together with 
\bel{Min8}
\left\langle \gamma
_s(\tau,\tau)+
\int_0^\tau\Psi(\sigma)\times
\bigl(\gamma
(\tau,\tau)-\gamma
(\tau,\sigma)\bigr)\, d\sigma+\bfv(\tau)\,,~\bfn(\tau,\tau)\right\rangle
~=~0\eeq
if $\gamma(\tau,\tau)\in \partial\Omega$.
By (\ref{RC}), $\bfv$ has the form 
\bel{vv2}\bfv(s)~=~\int_0^s \bar \omega(\sigma)\times 
\bigl(\gamma(t,s)-\gamma(t,\sigma)
\bigr))\, d\sigma\,,\eeq where the angular velocity is
\bel{o6}\bar\omega (\sigma)~=~-e^{-\beta(\tau-\sigma)} \int_{[\sigma,\tau]}  
\Big(\bfn(\tau,s')\times 
\bigl(\gamma(\tau,s')-\gamma(t,\sigma)\bigr)\Big)   d\mu(s')\,, \eeq
for some positive measure $\mu$ supported on the set $\chi(\tau)$. 
To achieve the proof we 
need to show that $\bar\omega (\cdot)$ provides the global minimizer 
for the optimization problem 
(\ref{Min1}) subject to the unilateral constraints (\ref{Min2})-(\ref{Min3}).
\v
{\bf 2.} Consider any other field of angular velocities, say 
$\bar\omega  +  \omega$.  The optimality of $\bar\omega $ will be proved by showing that 
\begi
\item~either $\E(\bar\omega +  \omega)\geq \E(\bar\omega )$,

\item or else, replacing $\bar\omega$ with $\bar\omega  +  \omega$, the 
constraints (\ref{Min2})-(\ref{Min3}) are no longer satisfied.
\endi
By the convexity of the integrand in (\ref{Min1}) it follows   
\bel{5}
{1\over 2}\int_0^\tau e^{\beta(\tau-\sigma)}|\bar\omega (\sigma)+ \omega(\sigma)|^2\, d\sigma~\geq~
{1\over 2}\int_0^\tau e^{\beta(\tau-\sigma)}|\bar\omega (\sigma)|^2\, d\sigma
+\int_0^\tau e^{\beta(\tau-\sigma)}\bigl\langle\bar\omega (\sigma),\, \omega(\sigma)
\bigr\rangle\, d\sigma\,.\eeq
The last term on the right hand side of (\ref{5}) is computed by
\bel{6}\bega{l}\ds
\int_0^\tau e^{\beta(\tau-\sigma)}\bigl\langle\bar\omega (\sigma),\, \omega(\sigma)
\bigr\rangle\, d\sigma\\[4mm]
\ds\qquad 
=~-\int_0^\tau 
e^{\beta(\tau-\sigma)}  \omega(\sigma)\cdot 
e^{-\beta(\tau-\sigma)}\left( \int_{[\sigma,\tau]}  
\Big(\bfn(\tau,s')\times 
\bigl(\gamma(\tau,s')-\gamma(t,\sigma)\bigr)\Big)   d\mu(s')\right)d\sigma\\[4mm]
\ds\qquad = ~\int_0^\tau\left(\int_{[\sigma,\tau]} 
\Big( \omega(\sigma)\times 
\bigl(\gamma(\tau,s')-\gamma(t,\sigma)\bigr) \Big)\cdot  \bfn(\tau,s')\,  d\mu(s')
\right) d\sigma.\enda\eeq

If $\bar\omega + \omega$ is admissible, then by (\ref{Min7})-(\ref{Min8})
it follows
$$\left(\int_0^{s'}
 \omega(\sigma)\times 
\bigl(\gamma(\tau,s')-\gamma(t,\sigma)\bigr)\,d\sigma\right)\cdot \bfn(\tau,s')
~\geq~0\qquad\forall s'\in \chi(\tau).$$
Integrating w.r.t.~$\mu$ and exchanging the order of integration one obtains
\bel{7}\bega{rl}0&\leq~\ds
\int_{[0,\tau]}\left(\int_0^{s'}
 \omega(\sigma)\times 
\bigl(\gamma(\tau,s')-\gamma(t,\sigma)\bigr)\,d\sigma\right) \cdot \bfn(\tau,s')\,d\mu(s')
\\[4mm]\ds
&=\ds ~\int_{0}^{\tau}\left(\int_{[\sigma,\tau]}  \omega(\sigma)\times 
\bigl(\gamma(\tau,s')-\gamma(t,\sigma)\bigr)\cdot \bfn(\tau,s') d\mu(s')\right)
\,d\sigma
\,.
\enda \eeq
Hence the right hand side of (\ref{6}) is nonnegative.

This shows that $\bar\omega (\cdot)$ in (\ref{o6}) provides the global minimizer to the
constrained optimization problem \eqref{Min1}-\eqref{Min3}. Since this minimization problem  has strictly convex cost and convex constraints, we conclude that 
$\bar\omega (\cdot)$ is the unique minimizer,  as  claimed in Theorem~3.
\endproof
\v

\end{document}